\documentclass[10pt,a4paper]{amsart}

\usepackage[marginpar=2cm,ignoremp,margin=3cm]{geometry}

\RequirePackage{doi}
\usepackage{hyperref}

\allowdisplaybreaks

\usepackage{amssymb}

\usepackage{braket}
\usepackage{mathtools}
\usepackage{thmtools}

\declaretheorem[
style=plain,
name=Theorem,
numberwithin=section,
refname={Theorem,Theorems},
Refname={Theorem,Theorems}
]{Thm}

\declaretheorem[
style=definition,
name=Example,
numberlike=Thm,
refname={Example,Examples},
Refname={Example,Examples}
]{Eg}
\declaretheorem[
style=plain,
name=Conjecture,
numberlike=Thm,
refname={Conjecture,Conjectures},
Refname={Conjecture,Conjectures}
]{Conj}
\declaretheorem[
style=plain,
name=Idea,
numberlike=Thm,
refname={Idea,Ideas},
Refname={Idea,Ideas}
]{Idea}
\declaretheorem[
style=definition,
name=Remark,
numberlike=Thm,
refname={Remark,Remarks},
Refname={Remark,Remarks}
]{Rem}

\newcommand{\Q}{\mathbb{Q}}
\newcommand{\C}{\mathbb{C}}

\newcommand{\mot}{\mathfrak{m}}

\DeclareMathOperator{\Ti}{Ti}
\newcommand{\abs}[1]{\left\lvert #1 \right\rvert}
\renewcommand{\Re}{\operatorname{Re}}

\newcommand{\poch}[2]{\left\{ #1 \right\}_{#2}}

\newmuskip\pFqmuskip
\newcommand*\pFq[6][8]{%
	\begingroup % only local assignments
	\pFqmuskip=#1mu\relax
	\mathchardef\normalcomma=\mathcode`,
	\mathcode`\,=\string"8000
	\begingroup\lccode`\~=`\,
	\lowercase{\endgroup\let~}\pFqcomma
	{}_{#2}F_{#3}{\left[\genfrac..{0pt}{}{#4}{#5};#6\right]}%
	\endgroup
}
\newcommand{\pFqcomma}{{\normalcomma}\mskip\pFqmuskip}

\renewcommand{\epsilon}{\varepsilon}
\newcommand{\sgnarg}[2]{
	\left( \genfrac{}{}{0pt}{0}{#1}{#2} \right)
}
\newcommand{\sgnargsm}[2]{
	\left( \textstyle\genfrac{}{}{0pt}{1}{#1}{#2} \right)
}

\newcommand{\ii}{\mathrm{i}}

\let\overlineO\overline
\renewcommand{\overline}[1]{\overlineO{\mathclap{\phantom{I}}#1}}

\makeatletter
\@namedef{subjclassname@2020}{%
	\textup{2020} Mathematics Subject Classification}
\makeatother

\begin{document}
	
	\title[Evaluation of \( t(\{\overline{1}\}^a, 1, \{\overline{1}\}^b) \)]{On the evaluation of the alternating  \\ multiple \( t \) value \( t(\{\overline{1}\}^a, 1, \{\overline{1}\}^b) \)}
	\author{Steven Charlton}
	\date{31 December 2021}
	
	
	\begin{abstract}
		We prove an evaluation for the stuffle-regularised multiple \( t \) value \( t^{\ast,V}(\{\overline{1}\}^a, 1, \{\overline{1}\}^b) \) in terms of \( \log(2), \zeta(k) \) and \( \beta(k) \).  This arises by evaluating the corresponding generating series using the Evans-Stanton/Ramanujan asymptotics of a zero-balanced hypergeometric function \( {}_3F_2 \), and an evaluation established by Li in an alternative approach to Zagier's evaluation of \( \zeta(\{2\}^a, 3, \{2\}^b) \).  We end with some discussion and conjectures on possible motivic applications.
	\end{abstract}
	
	\maketitle
	
	\section{Introduction and statement}
	
	The \emph{multiple zeta values} (MZV's) and \emph{multiple \( t \) values} (MtV's) with signs \( \epsilon_i \in S^1 = \{ z \in \C : \abs{z} = 1 \} \), and arguments \( k_i \in \mathbb{Z}_{>0} \), with \( k_d > 1 \), are defined by
	\begin{align*}
		\zeta\sgnarg{\epsilon_1,\ldots,\epsilon_d}{k_1,\ldots,k_d} &\coloneqq \sum_{0 < n_1 < \cdots < n_d} \frac{\epsilon_1^{n_1} \cdots \epsilon_d^{n_d}}{n_1^{k_1} \cdots n_d^{k_d}} \,, \\
		t\sgnarg{\epsilon_1,\ldots,\epsilon_d}{k_1,\ldots,k_d} &\coloneqq \sum_{0 < n_1 < \cdots < n_d} \frac{\epsilon_1^{n_1} \cdots \epsilon_d^{n_d}}{(2n_1-1)^{k_1} \cdots (2n_d-1)^{k_d}} \,.
	\end{align*}
	When all \( \epsilon_i \in \{\pm1\} \), write \( \overline{k_i} \) to denote arguments \( k_i \) which have associated sign \( \epsilon_i = -1 \).  Here \( d \) is called the \emph{depth} and \( k_1 + \cdots + k_d \) is called the \emph{weight} of the MZV or MtV.  The Dirichlet beta function is defined by 
	\[
	\beta(s) \coloneqq \sum_{n=0}^\infty \frac{(-1)^n}{(2n+1)^s} 
	\]
	so that \( t(\overline{n}) = -\beta(n) \), for \( n \in \mathbb{Z}_{>0} \).
	
	MZV's and by extension MtV's and other related objects are of significant interest in number theory (see \cite{zagier94}\cite{hoffman92} for the foundational results around MZV's although Euler already studied the case \( d \leq 2 \), and see \cite{hoffman19} for the recent (re-)introduction of MtV's following Nielsen's study of the case \( d = 1 \)).  They are also of interest for their applications to high energy physics (see \cite{broadhurstKreimer97} as a starting point).  Typically one is interested in understanding identities and relations between MZV's and MtV's, either particular cases or the general structures thereof.
	\medskip

	The main theorem of this note is an explicit evaluation and generating series for the stuffle-regularised \( t^{\ast,V}(\{\overline{1}\}^a,1,\{\overline{1}\}^b) \).  Therefore we briefly recall the idea of stuffle-regularisation of MZV's (see \cite{ikz06} for more details) and correspondingly MtV's.  The truncated \( \zeta_M(1) \) is well-known to satisfy
	\[
		\zeta_M(1) \coloneqq \sum_{n=1}^M \frac{1}{n} = \log(M) + \gamma + O\bigg(\frac{1}{M}\bigg) \,,
	\]
	where \( \gamma = 0.577\ldots \) is the Euler-Mascheroni constant.  So by application of the stuffle-product (for example \vspace{-0.5em}
	\begin{align*}
		\zeta(k_1) \zeta(\ell_1) &= \sum_{0 < n_1} \sum_{0 < m_1} \frac{1}{n_1^{k_1}} \cdot \frac{1}{m_1^{\ell_1}}  \\
		&=  \Bigg(\sum_{0 < n_1 <  m_1} + \sum_{0 < m_1  < n_1} + \sum_{0 < n_1 = m_1} \Bigg) \frac{1}{n_1^{k_1} m_1^{\ell_1}}  \\
		&= \zeta(k_1,\ell_1) + \zeta(\ell_1,k_1) + \zeta(\ell_1 + k_1) \,,
	\end{align*}
	in the case of two-single zeta values; when signs are present, they will be multiplied in the \( \ell_1 = k_1 \) term), one has by induction that any truncated MZV satisfies
	\[
		\zeta_M(k_1,\ldots,k_d) \coloneqq \sum_{0 < n_1 < \cdots < n_d < M} \frac{1}{n_1^{k_1} \cdots n_d^{k_d}} = Z^\ast(k_1,\ldots,k_d; \log(M) + \gamma) + O\bigg(\frac{\log^J{M}}{M}\bigg) \,,
	\]
	for some polynomial \( Z^\ast(k_1,\ldots,k_d; U) \) with convergent \( MZV \) coefficients, and some \( J \).  This polynomial defines the regularised version of \( \zeta^{\ast,U}(k_1,\ldots,k_d) \), with parameter \( \zeta^{\ast,U}(1) = U \); in particular \( \zeta^{\ast,U}(k_1,\ldots,k_d) = \zeta(k_1,\ldots,k_d) \), for \( k_d \neq 1 \).  Essentially, one can formally extend the stuffle-product to allow trailing 1's, then by considering, for \( k_d \neq 1 \), the difference
	\[
		\zeta^{\ast,U}(k_1,\ldots,k_d,\{1\}^\alpha) - \frac{1}{\alpha} \zeta^{\ast,U}(k_1,\ldots,k_d, \{1\}^{\alpha-1}) \zeta^{\ast,U}(1) \,,
	\]
	one obtains an expression with strictly fewer trailing 1's.  This writes every MZV as a polynomial in \( \zeta^{\ast,U}(1) \coloneqq U \) with convergent MZV coefficients.  The same process works for (alternating) MtV's, in particular, we have
	\begin{equation}
	\label{eqn:reg:tmm1}
	\begin{aligned}
		t^{\ast,V}(\{\overline{1}\}^m,1) &= t^{\ast,V}(\{\overline{1}\}^m) t^{\ast,V}(1) - \sum_{i=0}^{m-1} t^{\ast,V}(\{\overline{1}\}^{i}, 1, \{\overline{1}\}^{m-i})   - \sum_{i=0}^{m-1} t^{\ast,V}(\{\overline{1}\}^{i}, \overline{2}, \{\overline{1}\}^{m-1-i}) \\
		&= t(\{\overline{1}\}^m) V  - \sum_{i=0}^{m-1} t(\{\overline{1}\}^{i}, 1, \{\overline{1}\}^{m-i})   - \sum_{i=0}^{m-1} t(\{\overline{1}\}^{i}, \overline{2}, \{\overline{1}\}^{m-1-i}) \,.
	\end{aligned}
	\end{equation}

	The main Theorem of this note is now as follows.
	\begin{Thm}\label{thm:tmm1m}
		Consider the following generating series of the multiple \( t \) values \( t^{\ast,V}(\{\overline{1}\}^a, 1, \{\overline{1}\}^b) \), with stuffle-regularisation \( t^{\ast,V}(1) = V \) if necessary,
		\[
			F^V(x,y) = \sum_{a,b\geq0} (-1)^{a+b} t^{\ast,V}(\{\overline{1}\}^a,1,\{\overline{1}\}^b) x^a y^b \,.
		\]
		Then the generating series has the following closed form expression
		\begin{equation*}
		\begin{aligned}
		F^V(x,y) = & \frac{1}{2} \Big( \cos\Big(\frac{\pi x}{4}\Big) + \sin\Big(\frac{\pi x}{4}\Big) \Big)  \Big( A\Big(\frac{x-y}{4}\Big) - A\Big(\frac{x+y}{4}\Big) + 2 A\Big(\frac{x+y}{2}\Big) - \log(2) + 2 V \Big) \\
		& +  \frac{1}{2} \Big( \cos\Big(\frac{\pi y}{4}\Big) + \sin\Big(\frac{\pi y}{4}\Big) \Big) \Big( {-} A\Big(\frac{x-y}{8}\Big) + A\Big(\frac{x-y}{4}\Big) - 2  C\Big(\frac{x+y}{2}\Big) + \log(2) \Big) \,,
		\end{aligned}
		\end{equation*}
		where
		\begin{align*}
		A(z) &= \psi(1) - \frac{1}{2} (\psi(1+z) + \psi(1-z)) = \sum_{r=1}^\infty \zeta(2r+1) z^{2r} \\
		C(z) &= \frac{1}{8} \big( \psi(\tfrac{1}{4} + \tfrac{z}{4}) - \psi(\tfrac{1}{4} - \tfrac{z}{4}) - \psi(\tfrac{3}{4} + \tfrac{z}{4}) + \psi(\tfrac{3}{4} - \tfrac{z}{4}) \big) = \sum_{r=1}^\infty \beta(2r) z^{2r-1} \,,
		\end{align*}
		with \( \psi(x) = \frac{\mathrm{d}}{\mathrm{d}x} \log \Gamma(x) \) the digamma function.  In particular each \( t(\{\overline{1}\}^a,1,\{\overline{1}\}^b) \) is a polynomial in \( V \), Riemann zeta values \( \zeta(k) \), \( \log(2) \) and Dirichlet beta values \( \beta(k) \).
	\end{Thm}

	The following evaluation follows directly from the generating series, verifying the polynomial nature of the reduction
	\begin{align*}
	 t^{\ast,V}(\{\overline{1}\}^a,&1,\{\overline{1}\}^b) = {} \\[1ex]
	{} & \sum_{\substack{r = 2 \\ \text{even}}}^{a+b+1} \frac{(-1)^{\lfloor (a+b-r-1)/2 \rfloor}}{2(a+b-r)!} \Big(\frac{\pi}{4}\Big)^{a+b-r} \binom{r}{b} \left( \frac{1}{4^r} - \frac{(-1)^{b-r}}{4^r} - \frac{2}{2^r} \right) \zeta(r+1) \\
	& + \sum_{\substack{r = 2 \\ \text{even}}}^{a+b+1} \frac{(-1)^{\lfloor (a+b-r-1)/2 \rfloor}}{2(a+b-r)!} \Big(\frac{\pi}{4}\Big)^{a+b-r} \binom{r}{a} \left( \frac{(-1)^{a-r}}{8^r} - \frac{(-1)^{a-r}}{4^r} - \frac{2}{2^r} \right) \zeta(r+1) \\
	& - \sum_{\substack{r = 1 \\ \text{odd}}}^{a+b+1} \frac{(-1)^{\lfloor (a+b-r-1)/2 \rfloor}}{2(a+b-r)!} \Big(\frac{\pi}{4}\Big)^{a+b-r} \binom{r}{a} \left(  \frac{2}{2^r} \right) \beta(r+1) \\
	& + \delta_{a=0} \frac{(-1)^{\lfloor -(a+b)/2 \rfloor}}{2(a+b)!} \Big( \frac{\pi}{4} \Big)^{a+b} (2 V - \log(2))  + \delta_{b=0} \frac{(-1)^{\lfloor -(a+b)/2 \rfloor}}{2(a+b)!} \Big( \frac{\pi}{4} \Big)^{a+b} \log(2) \,.
	\end{align*}
	
		\subsection*{Acknowledgements} I am grateful to Michael Hoffman for frequent discussions on MtV's, MZV's and various related subjects during his extended stay at the Max-Planck Institute f\"ur Mathematik in Bonn, through spring and summer 2020.  I am also grateful to the MPIM for extended support, which facilitated these discussions with Michael Hoffman and precipitated the start of this work.  I am grateful to both Adam Keilthy and Danylo Radchenko for some discussions and suggestions on how to tackle the hypergeometric series which arise during the evaluation of the \( t(\{2\}^a,1,\{2\}^b) \) generating series which appears in \cite{charltont2212}, and directly informed the results of this note.  I was supported by DFG Eigene Stelle grant CH 2561/1-1, for Projektnummer 442093436.
	
	\section{Proof of \autoref{thm:tmm1m}}
	
	Firstly, recall that the \( {}_pF_{p-1} \) hypergeometric function is defined as
	\[
	\pFq{p}{p-1}{a_1,\ldots,a_p}{b_1,\ldots,b_{p-1}}{x} \coloneqq \sum_{m=0}^\infty \frac{\poch{a_1}{m} \cdots \poch{a_p}{m} }{\poch{b_1}{m} \cdots \poch{b_{p-1}}{m} } \frac{x^m}{m!} \,,
	\]
	where \( \poch{a}{m} = a(a+1) \cdots (a+m-1) \) is the ascending Pochhammer symbol.  Asymptotic and transformation properties of the \( {}_3F_2 \) function will play a key role in the proof of our generating series evaluation. \medskip
	
	Introduce the multivariable \( \Ti \)-function
	\[
		\Ti_{k_1,\ldots,k_d}(z_1,\ldots,z_d) \coloneqq \sum_{0 < n_1 < \cdots < n_d} \frac{z_1^{n_1} \cdots z_d^{n_d}}{(2n_1-1)^{k_1} \cdots (2n_d-1)^{k_d}} \,.
	\]
	(Note this is slightly different to the version in \cite{charltont2212}; here the arguments in the numerator have exponent \( n_i \), which is more suited to the case of alternating \( t \) values.)
	Then for \( (k_d,\epsilon_d) \neq (1,1) \), we find \( \Ti_{k_1,\ldots,k_d}(\epsilon_1,\ldots,\epsilon_d) \) evaluates to the alternating MtV with arguments \( k_1,\ldots, k_d \) and signs \( \epsilon_1,\ldots,\epsilon_d \).  We consider the generating series of \( \Ti_{\{1\}^{n+m+1}}(\{-1\}^m,z,\{-1\}^n) \) and seek to carefully evaluate the limit thereof as \( z \to 1^- \), somehow dealing with the divergence caused by \( \Ti_{1,\ldots,1}(-1,\ldots,-1,z) \xrightarrow{z\to1^-} t(\{\overline{1}\}^n,1) \).
	
	Consider
	\begin{align*}
		G(x,y;z) \coloneqq \sum_{a,b\geq0} (-1)^{a+b} \Ti_{\{1\}^{a+b+1}}(\{-1\}^a,z,\{-1\}^b) x^a y^b \\
				= \sum_{r=1}^\infty \prod_{0 < k < r} \Big( 1 - \frac{(-1)^k x}{2k-1} \Big) \frac{z^r}{2r-1} \prod_{\ell > r} \Big( 1 - \frac{(-1)^\ell x}{2\ell-1} \Big)
	\end{align*}
	Using standard evaluations one can show
	\[
		\prod_{k=1}^\infty \Big( 1 - \frac{(-1)^k y}{2k-1} \Big) = \cos\Big(\frac{\pi y}{4}\Big) + \sin\Big(\frac{\pi y}{4}\Big) \,,
	\]
	and so we can rewrite \( G(x,y;z) \) as follows.
	\[
		G(x,y;z) = \Big( \cos\Big(\frac{\pi y}{4}\Big) + \sin\Big(\frac{\pi y}{4}\Big) \Big) \cdot \sum_{r=1}^\infty \prod_{0 < k < r} \Big( 1 - \frac{(-1)^k x}{2k-1} \Big) \frac{z^r}{2r-1} \prod_{0 < \ell \leq r} \Big( 1 - \frac{(-1)^\ell x}{2\ell-1} \Big)^{-1}
	\]
	We split the summation in \( G(x,y;z) \) into odd- and even- indexed terms, and sum each separately.  We can check the following (where \( \poch{x}{k} = x (x+1) \cdots (x+k-1) \) is, as already mentioned, the ascending-Pochhammer symbol): 
	\begin{align*}
		 (r = 2m + 2) \quad \quad & \prod_{0 < k < r} \Big( 1 - \frac{(-1)^k x}{2k-1} \Big) \frac{z^r}{2r-1} \prod_{0 < \ell \leq r} \Big( 1 - \frac{(-1)^\ell x}{2\ell-1} \Big)^{-1} \\
		 & {} = \frac{(1+x)z^2}{(3-y)(1+y)} \cdot \frac{
		 	\poch{1}{m}
		 	\poch{\tfrac{3}{4} - \tfrac{x}{4}}{m}
		 	\poch{\tfrac{5}{4} + \tfrac{x}{4}}{m}
	 	}{
	 		\poch{\tfrac{7}{4} - \tfrac{y}{4}}{m}
	 		\poch{\tfrac{5}{4} + \tfrac{y}{4}}{m}
	 	} \frac{(z^2)^m}{m!} \,, \\
 	 (r = 2m + 1) \quad \quad & \prod_{0 < k < r} \Big( 1 - \frac{(-1)^k x}{2k-1} \Big) \frac{z^r}{2r-1} \prod_{0 < \ell \leq r} \Big( 1 - \frac{(-1)^\ell x}{2\ell-1} \Big)^{-1} \\
 	& {} = \frac{z}{1+y} \cdot \frac{
 		\poch{1}{m}
 		\poch{\tfrac{3}{4} - \tfrac{x}{4}}{m}
 		\poch{\tfrac{1}{4} + \tfrac{x}{4}}{m}
 	}{
 		\poch{\tfrac{3}{4} - \tfrac{y}{4}}{m}
 		\poch{\tfrac{5}{4} + \tfrac{y}{4}}{m}
 	} \frac{(z^2)^m}{m!} \,.
	\end{align*}
	Both of these are in exactly the right form to sum to a \( _{3}F_{2} \)-hypergeometric function (the summation index in both cases starts at \( m = 0 \)), and so we obtain
	\begin{align*}
		G(x,y;z) = \Big( \cos\Big(\frac{\pi y}{4}\Big) + \sin\Big(\frac{\pi y}{4}\Big) \Big) \cdot \Bigg\{
			& \frac{(1+x)z^2}{(3-y)(1+y)} \cdot \pFq{3}{2}{1,\tfrac{3}{4} - \tfrac{x}{4},\tfrac{5}{4} + \tfrac{x}{4}}{
				\tfrac{7}{4} - \tfrac{y}{4},\tfrac{5}{4} + \tfrac{y}{4}}{z^2} 
			\\
			& {} + \frac{z}{1+y} \cdot \pFq{3}{2}{1,\tfrac{3}{4} - \tfrac{x}{4},\tfrac{1}{4} + \tfrac{x}{4}}{
				\tfrac{3}{4} - \tfrac{y}{4},\tfrac{5}{4} + \tfrac{y}{4}}{z^2}
		\Bigg\} \,.
	\end{align*}
	We now consider how the divergence in \( \Ti_{\{1\}^{m+1}}(\{-1\}^m, z) \) arises, so we can compensate for it in the generating series.  We have by the stuffle-multiplication of \( \Ti \), that
	\begin{align*}
		\Ti_{\{1\}^m,1}(\{-1\}^m, z) = & {} \Ti_{\{1\}^m}(\{-1\}^m) \Ti_1(z) \\
			& - \sum_{j=0}^{m-1} \Ti_{\{1\}^{m+1}}(\{-1\}^a, z, \{1\}^{m-a}) \\
			& - \sum_{j=0}^{m-1} \Ti_{\{1\}^a,2,\{1\}^{m-1-a}}(\{-1\}^a, -z, \{-1\}^{m-1-a}) \,.
	\end{align*}
	On rearranging and taking the generating series of both sides, we find that
	\begin{align*}
	& \sum_{m\geq0} (-1)^m \Ti_{\{1\}^m,1}(\{-1\}^m, z) x^m - \Ti_1(z) \sum_{m\geq0} (-1)^m \Ti_{\{1\}^m}(\{-1\}^m) x^m \\
	& = -\sum_{m\geq0} (-1)^m \Big( \sum_{j=0}^{m-1} \Ti_{\{1\}^{m+1}}(\{-1\}^j, z, \{1\}^{m-j}) \\ & \quad\quad\quad\quad + \sum_{j=0}^{m-1} \Ti_{\{1\}^j,2,\{1\}^{m-1-j}}(\{-1\}^j, -z, \{-1\}^{m-1-j}) \Big) x^{m} \\
	& \xrightarrow{z\to1^-} -\sum_{m\geq0} (-1)^m \Big( \sum_{j=0}^{m-1} t(\{\overline{1}\}^j, 1, \{\overline{1}\}^{m-j}) \\ & \quad\quad\quad\quad + \sum_{j=0}^{m-1} t(\{\overline{1}\}^j,\overline{2},\{\overline{1}\}^{m-1-j}) \Big) x^{m} \\ 
	&= \sum_{m\geq0} (-1)^m t^{\ast,V=0}(\{\overline{1}\}^m, 1) x^m \,,
	\end{align*}
	as \( z \to 1^{-} \).  Here we have recognized the limiting generating series as that of \( t^{\ast,V=0}(\{\overline{1}\}^m, 1) \), where \( t^{\ast,V=0} \) represents the stuffle-regularisation polynomial of \( t(\{\overline{1}\}^m,1) \) evaluated at \( V = 0 \).  (This is the regularisation where we take \( t^{\ast,V=0}(1) = 0 \), cf. \eqref{eqn:reg:tmm1}.)  We also notice that
	\[
		\Ti_1(z) = \sqrt{z} \tanh^{-1}(\sqrt{z}) \,,
	\]
	and from \cite[Corollary 6.1]{hoffman19}, the generating series of \( t(\{\overline{1}\}^m) \) is given by
	\begin{equation}\label{eqn:tmmm_gs}
		\sum_{m\geq0} (-1)^m t(\{\overline{1}\}^m) x^m = \cos\Big(\frac{\pi x}{4}\Big) + \sin\Big(\frac{\pi x}{4}\Big) \,.
	\end{equation}
	We hence find the limit of the following combination gives the generating series of \( t^{\ast,V=0}(\{\overline{1}\}^a, 1, \{\overline{1}\}^b) \), with a suitably stuffle-regularised variant \( t^{\ast,V=0}(\{\overline{1}\}^a,1) \) in the case \( b = 0 \).
	\begin{align*}& \lim_{z\to1^-} G(x,y;z) - \sqrt{z} \tanh^{-1}(\sqrt{z}) \cdot \Big( \cos\Big(\frac{\pi x}{4}\Big) + \sin\Big(\frac{\pi x}{4}\Big) \Big) \\
		& = \sum_{a\geq0,b>0} (-1)^{a+b} t(\{\overline{1}\}^a,1,\{\overline{1}\}^b) x^a y^b + \sum_{a \geq 0} (-1)^a t^{\ast,V=0}(\{\overline{1}\}^a, 1) x^a \\
		& \eqqcolon F(x,y) \,.
	\end{align*}
	
	We now want to take the limit \( \lim_{z\to1^-} G(x,y;z) - \sqrt{z} \tanh^{-1}(\sqrt{z}) \cdot \big( \cos\big(\frac{\pi x}{4}\big) + \sin\big(\frac{\pi x}{4}\big) \big) \).  For this, we recall the following result \cite[Theorem 3]{evans84} (proving a claim of Ramanujan), which treats the asymptotics of the 0-balanced hypergeometric function \( {}_{3}F_{2} \), which both of the hypergeometric series in \( G(x,y;z) \) are.
	
	\begin{Thm}[Evans-Stanton 1984 \cite{evans84}, Ramanujan]
		If \( a + b + c = d + e \), and \( \Re(c) > 0 \), then as \( u \to 1^- \),
		\[
		\frac{\Gamma(a)\Gamma(b)\Gamma(c)}{\Gamma(d)\Gamma(e)} \cdot \pFq{3}{2}{a,b,c}{d,e}{u} = -\log(1-u) + L + O((1-u)\log(1-u)) \,,
		\]
		where
		\[
		L = -2\gamma - \psi(a) - \psi(b) + \sum_{k=1}^\infty \frac{\poch{d-c}{k} \poch{e-c}{k}}{\poch{a}{k} \poch{b}{k} k} \,.
		\]
		Here \( \gamma \approx 0.577\ldots \) is the Euler-Mascheroni constant, \( \psi(x) = \frac{\mathrm{d}}{\mathrm{d}x} \log \Gamma(x) \) is the digamma function, and \( \poch{x}{k} = x(x+1) \ldots (x + k-1) \) is the ascending Pochhammer symbol.
	\end{Thm}
	
	If we apply it to both of the above \( {}_3F_{2} \) functions, with \( c = 1 \) after appropriately permuting the arguments, we find (after some simplification of the prefactor and the resulting gamma function combination) that
	\begin{align*}
		& \Big( \cos\Big(\frac{\pi y}{4}\Big) + \sin\Big(\frac{\pi y}{4}\Big) \Big) \cdot \frac{(1+x)z^2}{(3-y)(1+y)} \cdot \pFq{3}{2}{1,\tfrac{3}{4} - \tfrac{x}{4},\tfrac{5}{4} + \tfrac{x}{4}}{
			\tfrac{7}{4} - \tfrac{y}{4},\tfrac{5}{4} + \tfrac{y}{4}}{z^2} = \\
		 & \frac{1}{4} z^2 \Big( \cos\Big(\frac{\pi x}{4}\Big) + \sin\Big(\frac{\pi x}{4}\Big) \Big) \bigg(
		 	-2\gamma - \psi(\tfrac{3}{4} - \tfrac{x}{4})  - \psi(\tfrac{5}{4} + \tfrac{x}{4}) - \log(1-z^2) \\
		 & \hspace{10em} + \sum_{k=1}^\infty \frac{
		 		\poch{\tfrac{3}{4}-\tfrac{y}{4}}{k}
	 			\poch{\tfrac{1}{4}+\tfrac{y}{4}}{k}
 		}{k 
 			\poch{\tfrac{3}{4}-\tfrac{x}{4}}{k}
 			\poch{\tfrac{5}{4}+\tfrac{x}{4}}{k}
 			} \bigg) + O((1 - z)\log(1 - z))
 		\shortintertext{and}
 		& \Big( \cos\Big(\frac{\pi y}{4}\Big) + \sin\Big(\frac{\pi y}{4}\Big) \Big) \cdot \frac{z}{1+y} \cdot \pFq{3}{2}{1,\tfrac{3}{4} - \tfrac{x}{4},\tfrac{1}{4} + \tfrac{x}{4}}{
 			\tfrac{3}{4} - \tfrac{y}{4},\tfrac{5}{4} + \tfrac{y}{4}}{z^2} = \\
 		& \frac{1}{4} z \Big( \cos\Big(\frac{\pi x}{4}\Big) + \sin\Big(\frac{\pi x}{4}\Big) \Big) \bigg(
 		-2\gamma - \psi(\tfrac{3}{4} - \tfrac{x}{4})  - \psi(\tfrac{1}{4} + \tfrac{x}{4}) - \log(1-z^2) \\
 		& \hspace{10em} + \sum_{k=1}^\infty \frac{
 			\poch{-\tfrac{1}{4}-\tfrac{y}{4}}{k}
 			\poch{\tfrac{1}{4}+\tfrac{y}{4}}{k}
 		}{k 
 			\poch{\tfrac{3}{4}-\tfrac{x}{4}}{k}
 			\poch{\tfrac{1}{4}+\tfrac{x}{4}}{k}
 		} \bigg) + O((1 - z)\log(1 - z))
	\end{align*}
	Since
	\begin{align*}
	& \lim_{z\to1^-} \Big( \cos\Big(\frac{\pi x}{4}\Big) + \sin\Big(\frac{\pi x}{4}\Big) \Big) \cdot \Big( \frac{1}{4} z^2 \log(1 - z^2) + \frac{1}{4} z \log(1-z^2) + \sqrt{z} \tanh^{-1}(\sqrt{z}) \Big)  \\
	& \hspace{5em} {} = \frac{3}{2} \log(2) \Big( \cos\Big(\frac{\pi x}{4}\Big) + \sin\Big(\frac{\pi x}{4}\Big) \Big) \,,
	\end{align*}
	we find these asymptotic formulae lead to following limit for \( G(x,y;z) \) that we seek
	\begin{align*}
		F(x,y) = {} & \lim_{z\to1^-} G(x,y;z) - \sqrt{z} \tanh^{-1}(\sqrt{z}) \cdot \Big( \cos\Big(\frac{\pi x}{4}\Big) + \sin\Big(\frac{\pi x}{4}\Big) \Big) \\
		= {} & \frac{1}{4}  \Big( \cos\Big(\frac{\pi x}{4}\Big) + \sin\Big(\frac{\pi x}{4}\Big) \Big) \bigg( -4 \gamma - 6 \log(2)  - \psi(\tfrac{1}{4} + \tfrac{x}{4}) - 2\psi(\tfrac{3}{4} - \tfrac{x}{4}) - \psi(\tfrac{5}{4} + \tfrac{x}{4}) \\
		& \hspace{5em} + \sum_{k=1}^\infty \frac{
			\poch{\tfrac{3}{4}-\tfrac{y}{4}}{k}
			\poch{\tfrac{1}{4}+\tfrac{y}{4}}{k}
		}{k 
			\poch{\tfrac{3}{4}-\tfrac{x}{4}}{k}
			\poch{\tfrac{5}{4}+\tfrac{x}{4}}{k}
		} +  \sum_{k=1}^\infty \frac{
			\poch{-\tfrac{1}{4}-\tfrac{y}{4}}{k}
			\poch{\tfrac{1}{4}+\tfrac{y}{4}}{k}
		}{k 
			\poch{\tfrac{3}{4}-\tfrac{x}{4}}{k}
			\poch{\tfrac{1}{4}+\tfrac{x}{4}}{k}
		}
	\bigg)
	\end{align*}
	Using the same observation as in Proposition 1 of \cite{zagier2232}, we may rewrite these Pochhammer sums in a more useful way.  Specifically we have
	\[
		\sum_{k=1}^\infty \frac{\poch{a}{k}\poch{b}{k} }{k \poch{c}{k}\poch{d}{k}} = \frac{\mathrm{d}}{\mathrm{d}Z} \bigg|_{Z=0} \pFq{3}{2}{a,b,Z}{c,d}{1} \,,
	\]
	so it follows that
	\begin{align*}
		F(x,y) =  & \frac{1}{4}  \Big( \cos\Big(\frac{\pi x}{4}\Big) + \sin\Big(\frac{\pi x}{4}\Big) \Big) \Bigg\{ {-}4 \gamma - 6 \log(2)  - \psi(\tfrac{1}{4} + \tfrac{x}{4}) - 2\psi(\tfrac{3}{4} - \tfrac{x}{4}) - \psi(\tfrac{5}{4} + \tfrac{x}{4}) \\
		& \hspace{5em} + \frac{\mathrm{d}}{\mathrm{d}z} \bigg\rvert_{z=0} \Bigg( \pFq{3}{2}{
			\tfrac{3}{4}-\tfrac{y}{4},\tfrac{1}{4}+\tfrac{y}{4},z
			}{
			\tfrac{3}{4}-\tfrac{x}{4},\tfrac{5}{4}+\tfrac{x}{4}
			}{1}
		+ \pFq{3}{2}{
			-\tfrac{1}{4}-\tfrac{y}{4},\tfrac{1}{4}+\tfrac{y}{4},z
		}{
			\tfrac{3}{4}-\tfrac{x}{4},\tfrac{1}{4}+\tfrac{x}{4}
		}{1}	
		\Bigg) \! \Bigg\} \,.
	\end{align*}
		Apply to the second term following contiguous function relation (with arguments in this order)
	\[
	(a-b)p \cdot \pFq{3}{2}{a,b,c}{p,q}{z} - b(a-p) \cdot \pFq{3}{2}{a,b+1,c}{p+1,q}{z} + a(b-p) \cdot \pFq{3}{2}{a+1,b,c}{p+1,q}{z} = 0 \,,
	\] and we find
	\begin{equation}\label{eqn:fxy_asder}
	\begin{aligned}
		& F(x,y) = \frac{1}{4}  \Big( \cos\Big(\frac{\pi x}{4}\Big) + \sin\Big(\frac{\pi x}{4}\Big) \Big) \Bigg\{ {-}4 \gamma - 6 \log(2)  - \psi(\tfrac{1}{4} + \tfrac{x}{4}) - 2\psi(\tfrac{3}{4} - \tfrac{x}{4}) - \psi(\tfrac{5}{4} + \tfrac{x}{4}) \\
		& \hspace{5em} + \frac{\mathrm{d}}{\mathrm{d}Z} \bigg\rvert_{Z=0} \Bigg( \pFq{3}{2}{
			\tfrac{3}{4}-\tfrac{y}{4},\tfrac{1}{4}+\tfrac{y}{4},Z
		}{
			\tfrac{3}{4}-\tfrac{x}{4},\tfrac{5}{4}+\tfrac{x}{4}
		}{1}
		+ \frac{x-y-4}{2(x-3)}\cdot\pFq{3}{2}{
			-\tfrac{1}{4}-\tfrac{y}{4},\tfrac{5}{4}+\tfrac{y}{4},Z
		}{
			\tfrac{7}{4}-\tfrac{x}{4},\tfrac{1}{4}+\tfrac{x}{4}
		}{1} \\
	& \hspace{10.5em} + \frac{x+y-2}{2(x-3)}\cdot\pFq{3}{2}{
			\tfrac{3}{4}-\tfrac{y}{4},\tfrac{1}{4}+\tfrac{y}{4},Z
			}{
			\tfrac{7}{4}-\tfrac{x}{4},\tfrac{1}{4}+\tfrac{x}{4}
			}{1}	
		\Bigg) \! \Bigg\} \,.
	\end{aligned}
	\end{equation}
	All three terms are of the following form, with \( X = \tfrac{1}{4} + \tfrac{y}{4}, Y = \frac{1}{4} + \frac{x}{4} \), then \( X = -\frac{1}{4}-\tfrac{y}{4} , Y = \frac{3}{4} - \tfrac{x}{4} \) and then \( X = \frac{1}{4}+\tfrac{y}{4} , Y = \frac{3}{4} - \tfrac{x}{4} \) respectively:
	\begin{equation}\label{eqn:f32_eval}
	\begin{aligned}
	& \frac{\mathrm{d}}{\mathrm{d}Z} \Bigg|_{Z=0} \pFq{3}{2}{X,1-X,Z}{1-Y,1+Y}{1} \\
	& {} = \psi(1 + Y) + \psi(1 - Y) - \psi(1 - X + Y) - \psi(1 - X - Y) \\
	& \quad\quad - \frac{\sin(\pi X)}{\sin(\pi Y)} \cdot \Big[ \psi(1 - X + Y) - \psi(1 - X - Y) - \psi\big( 1 - \tfrac{X-Y}{2} \big) + \psi\big( 1 - \tfrac{X+Y}{2} \big) \Big] \,.
	\end{aligned}
	\end{equation}
	This evaluation is given in Equation 9 of \cite{li13}, and hence reduces \( F(x,y) \) to a combination of digamma functions.
	
	After substituting in the indicated \( X, Y \) values, and some amount of simplification (see \autoref{rem:fxy_simp} below), we find
	\begin{equation}\label{eqn:fxy_as_ac}
	\begin{aligned}
		F(x,y) = & \frac{1}{2} \Big( \cos\Big(\frac{\pi x}{4}\Big) + \sin\Big(\frac{\pi x}{4}\Big) \Big) \Big( A\Big(\frac{x-y}{4}\Big) - A\Big(\frac{x+y}{4}\Big) + 2 A\Big(\frac{x+y}{2}\Big) - \log(2) \Big) \\
		& +  \frac{1}{2} \Big( \cos\Big(\frac{\pi y}{4}\Big) + \sin\Big(\frac{\pi y}{4}\Big) \Big) \Big( {-} A\Big(\frac{x-y}{8}\Big) + A\Big(\frac{x-y}{4}\Big) - 2  C\Big(\frac{x+y}{2}\Big) + \log(2) \Big) \,,
	\end{aligned}
	\end{equation}
where \( A \) and \( C \) are certain generating series of \( \zeta \)- and \( \beta \)-values as defined below.  The series \( A \) is defined as in \cite{zagier2232}, namely
	\begin{align*}
		A(z) &\coloneqq \psi(1) - \frac{1}{2} (\psi(1+z) + \psi(1-z)) = \sum_{r=1}^\infty \zeta(2r+1) z^{2r} \,.
	\end{align*}
	The generating \( C \) is defined as follows, to give an analogous Dirichlet-\( \beta \) generating series:
	\[
		C(z) \coloneqq \frac{1}{8} \big( \psi(\tfrac{1}{4} + \tfrac{z}{4}) - \psi(\tfrac{1}{4} - \tfrac{z}{4}) - \psi(\tfrac{3}{4} + \tfrac{z}{4}) + \psi(\tfrac{3}{4} - \tfrac{z}{4}) \big) = \sum_{r=1}^\infty \beta(2r) z^{2r-1} \,.
	\]
	This formula for \( C \) follows by interchanging the summation, and applying partial fractions to the result, in 
	\begin{align*}
		\sum_{r=1}^\infty \beta(2r)z^{2r-1} &= \sum_{r=1}^\infty \sum_{k=1}^\infty \frac{-(-1)^r}{(2k-1)^{2r}} x^{2r-1} \\
		&= \sum_{m=0}^\infty \frac{(-1)^m}{2(1-2m+x)} - \frac{(-1)^m}{2(1-2m-x)} \,.
	\end{align*}
	
	\begin{Rem}\label{rem:fxy_simp}
		Since the steps of simplification of \( F(x,y) \) are rather involved, we indicate some strategy to follow in the verification process.  The goal is to check \( I = \eqref{eqn:fxy_asder} - \eqref{eqn:fxy_as_ac} = 0 \), where \eqref{eqn:fxy_asder} is evaluated via the formula for \( \frac{\mathrm{d}}{\mathrm{d}z} \) of a certain \( {}_3F_2 \) hypergeometric series given in \eqref{eqn:f32_eval}. 
		
		Via the functional equation
		\[
			\psi(z + 1) - \psi(z) = \frac{1}{z} \,,
		\]
		every argument of \( (x-3)I \) (we multiply by $(x-3)$ to eliminate rational function coefficients at this point) can be reduced to one of the form \( \alpha + \beta x + \gamma y \), where \( \alpha = -\tfrac{1}{4}, 0, \frac{1}{4}, \frac{1}{2} \).  Using the functional equation
		\begin{equation}\label{eqn:psi:mirror}
			\psi(-z) = \psi(z) + \frac{1}{z} + \pi \cot(\pi z) \,,
		\end{equation}
		we can also reduce the \( \alpha = -\tfrac{1}{4} \) arguments to \( \tfrac{1}{4} \).	In particular we are left with only the following arguments
		\begin{align*}
			\big\{ &
			\psi \left(-\tfrac{x}{2}-\tfrac{y}{2}\right),
			\psi \left(-\tfrac{x}{4}-\tfrac{y}{4}\right),
			\psi \left(\tfrac{x}{4}-\tfrac{y}{4}\right),
			\psi\left(\tfrac{1}{4}-\tfrac{x}{8}-\tfrac{y}{8}\right),
			\psi \left(\tfrac{x}{8}-\tfrac{y}{8}\right),
			\psi \left(-\tfrac{x}{8}+\tfrac{y}{8}\right),\\
			& \psi \left(\tfrac{1}{2}-\tfrac{x}{8}+\tfrac{y}{8}\right),
			\psi\left(\tfrac{1}{4}+\tfrac{x}{8}+\tfrac{y}{8}\right),
			\psi\left(-\tfrac{x}{4}+\tfrac{y}{4}\right),
			\psi\left(\tfrac{x}{4}+\tfrac{y}{4}\right),
			\psi\left(\tfrac{1}{2}+\tfrac{x}{4}+\tfrac{y}{4}\right),
			\psi\left(\tfrac{x}{2}+\tfrac{y}{2}\right) \big\}
		\end{align*}
		Applying the symmetry in \eqref{eqn:psi:mirror} further times, reduces \( -\tfrac{x}{2}-\tfrac{y}{2}, -\tfrac{x}{4}-\tfrac{y}{4}, -\tfrac{x}{4}+\tfrac{y}{4}, -\tfrac{x}{8}+\tfrac{y}{8} \) to their positive counterparts.
		Using \eqref{eqn:psi:mirror} (again, undoing the previous step in one case!) and the functional equation
		\begin{equation}\label{eqn:psi:dup}
			\psi(2z) = \tfrac{1}{2} \big( \psi(z + \tfrac{1}{2}) + \psi(z) \big) + \log(2) \,,
		\end{equation}
		with
		\( z = z' - \frac{1}{4} \), we can express \( \psi(\tfrac{1}{4} - \tfrac{x}{8} - \tfrac{y}{8}) \) in terms of \( \psi(-\tfrac{1}{4} + \tfrac{x}{8} + \tfrac{y}{8}) \) first, then in terms of \( \psi(\tfrac{1}{4} + \tfrac{x}{8} + \tfrac{y}{8}) \) and \( \psi(-\tfrac{1}{2} - \tfrac{x}{4} - \tfrac{y}{4}) \).
		
		Use \eqref{eqn:psi:mirror} to reduce \( \psi(-\tfrac{1}{2} + z) \) to \( \psi(\tfrac{1}{2} + z) \), and apply \eqref{eqn:psi:dup} to reduce this to \( \psi(z) \) and \( \psi(2z) \), where necessary.  Then using \eqref{eqn:psi:mirror} we can finally reduce all arguments to the following set of 5 possibilities
		\[
			\left\{
			\psi \left(\tfrac{x}{4}-\tfrac{y}{4}\right),
			\psi \left(\tfrac{x}{8}-\tfrac{y}{8}\right),
			\psi\left(\tfrac{1}{4} +\tfrac{x}{8}+\tfrac{y}{8}\right),
			\psi \left(\tfrac{x}{4}+\tfrac{y}{4}\right),
			\psi \left(\tfrac{x}{2}+\tfrac{y}{2}\right)
			\right\} \,.
		\]
		At this point the coefficients of each of these 5 arguments is a trigonometric rational function, as is the constant coefficient of the \( \psi \)-polynomial.  Each of these trig-functions can be written as a rational function (with \(\sqrt{2}\)-coefficients) in 
		\[
			\cos\Big(\frac{\pi x}{8}\Big),			\sin\Big(\frac{\pi x}{8}\Big), 
			\cos\Big(\frac{\pi y}{8}\Big),			\sin\Big(\frac{\pi y}{8}\Big) \,,
		\]
		using the addition formulae for \( \sin,\cos \), and then reducing \( \frac{\pi x}{2} \)- and \( \frac{\pi x}{4} \)-arguments to the above via the double-angle formulae.  (For cosine, \( \cos(2x) = \cos(x)^2 - \sin(x)^2 \) seems to be the more useful variant.)  At this point the coefficients of \( \psi \)'s vanish identically (without relating powers of \( \sin \) and \( \cos \)), whereas the constant coefficient factors into an expression involving \( 
		\cos^2(\frac{\pi x}{8}) + \sin^2(\frac{\pi x}{8}) - \cos^2(\frac{\pi y}{8}) - \sin^2(\frac{\pi y}{8}) \), which is \( 0 \) via Pythagoras.  This shows that \( \eqref{eqn:fxy_asder} - \eqref{eqn:fxy_as_ac} = 0 \) as claimed, hence verifying the expression for \( F(x,y) \) given in \eqref{eqn:fxy_as_ac}.
	\end{Rem}

	At this point we have the generating series involving (for \( b = 0 \)) the regularised values \( t^{\ast,V=0} \) at \( V=0 \), however for M$t$V's the most natural stuffle-regularisation would seem to be \( V = \log(2) \).  Fortunately, since \( t(\{\overline{1}\}^m,1) \) only involves a single argument 1, hence the regularisation polynomial is given (cf. \eqref{eqn:reg:tmm1}) by
	\[
		t^{\ast,V}(\{\overline{1}\}^m,1) = t(\{\overline{1}\}^m) V + t^{\ast,V=0}(\{\overline{1}\}^m, 1) \,
	\]
	meaning we can (re-)construct the entire polynomial from this constant term.  In particular (see \eqref{eqn:tmmm_gs})
	\[
	\sum_{m\geq0} (-1)^m t^{\ast,V}(\{\overline{1}\}^m,1) x^m = V \Big( \cos\Big(\frac{\pi x}{4}\Big)+ \sin\Big(\frac{\pi x}{4}\Big)\Big)   + \sum_{m\geq0} (-1)^m t^{\ast,V=0}(\{\overline{1}\}^m, 1) x^m \,
	\]
	
	We can therefore give the arbitrary stuffle-regularisation generating series as follows
	\begin{equation}\label{eqn:fxygen_as_ac}
	\begin{aligned}
	F^V(x,y) = & \frac{1}{2} \Big( \cos\Big(\frac{\pi x}{4}\Big) + \sin\Big(\frac{\pi x}{4}\Big) \Big)  \Big( A\Big(\frac{x-y}{4}\Big) - A\Big(\frac{x+y}{4}\Big) + 2 A\Big(\frac{x+y}{2}\Big) - \log(2) + 2 V \Big) \\
	& +  \frac{1}{2} \Big( \cos\Big(\frac{\pi y}{4}\Big) + \sin\Big(\frac{\pi y}{4}\Big) \Big) \Big( {-} A\Big(\frac{x-y}{8}\Big) + A\Big(\frac{x-y}{4}\Big) - 2  C\Big(\frac{x+y}{2}\Big) + \log(2) \Big) \,,
	\end{aligned}
	\end{equation}
	
	The final step is to extract an explicit evaluation for \( t(\{\overline{1}\}^a,1,\{\overline{1}\}^b) \) from this generating series.  One can check easily that
	\[
		[x^a y^b] \sum_{i=0}^\infty f(i) x^i \cdot \sum_{j=0}^\infty g(j) (x+y)^j = \sum_{n=0}^{a+b} \binom{n}{b} g(a+b-n) f(n) \,,
	\]
	where \( [x^a y^b] \) denotes the coefficient of \( x^a y^a \) in the terms thereafter.  So we can readily extract the following formulae from \( F^V(x,y) \).
	\begin{align*}
		t^{\ast,V}(\{\overline{1}\}^a,&1,\{\overline{1}\}^b) = {} \\[1ex] & \sum_{\substack{r = 2 \\ \text{even}}}^{a+b+1} \frac{(-1)^{\lfloor (a+b-r-1)/2 \rfloor}}{2(a+b-r)!} \Big(\frac{\pi}{4}\Big)^{a+b-r} \binom{r}{b} \left( \frac{1}{4^r} - \frac{(-1)^{b-r}}{4^r} - \frac{2}{2^r} \right) \zeta(r+1) \\
		 & + \sum_{\substack{r = 2 \\ \text{even}}}^{a+b+1} \frac{(-1)^{\lfloor (a+b-r-1)/2 \rfloor}}{2(a+b-r)!} \Big(\frac{\pi}{4}\Big)^{a+b-r} \binom{r}{a} \left( \frac{(-1)^{a-r}}{8^r} - \frac{(-1)^{a-r}}{4^r} - \frac{2}{2^r} \right) \zeta(r+1) \\
		 & - \sum_{\substack{r = 1 \\ \text{odd}}}^{a+b+1} \frac{(-1)^{\lfloor (a+b-r-1)/2 \rfloor}}{2(a+b-r)!} \Big(\frac{\pi}{4}\Big)^{a+b-r} \binom{r}{a} \left(  \frac{2}{2^r} \right) \beta(r+1) \\
		 & + \delta_{a=0} \frac{(-1)^{\lfloor -(a+b)/2 \rfloor}}{2(a+b)!} \Big( \frac{\pi}{4} \Big)^{a+b} (2 V - \log(2))  + \delta_{b=0} \frac{(-1)^{\lfloor -(a+b)/2 \rfloor}}{2(a+b)!} \Big( \frac{\pi}{4} \Big)^{a+b} \log(2) \,.
	\end{align*}
	This completes the proof of \autoref{thm:tmm1m}, and the explicit evaluation thereafter. \hfill \qedsymbol	

	\section{Potential motivic applications}

	In \cite{brown12} and \cite{zagier2232} with the case \( \zeta(\{2\}^a, 3, \{2\}^b) \), in \cite{murakami21} with \( t(\{2\}^a,3,\{2\}^b) \) and in \cite{charltont2212} with the case \( t^{\ast,V=0}(\{2\}^a, 1, \{2\}^b) \), the various identities (and the arithmetic of their coefficients) were applied to show linear independence and/or basis results on the motivic level.  It should therefore be possible to lift \autoref{thm:tmm1m} to a motivic version, and show on the motivic level some independence and/or basis results for alternating MtV's.
	
	Equivalently (after extending coefficients to \( \mathbb{Q}(\ii) \)), one should also obtain results about coloured MZV's of level \( N = 4 \) (i.e. \( \epsilon_i \in \{ \pm 1, \pm \ii \} \), roots of unity of order \( N = 4 \)).  In general we have the following expression for MtV's of level \( N \) in terms of coloured MZV's of level \( 2N  \) (after fixing some choice of square roots, which we symmetrise over anyway)
	\begin{align*}
		t\sgnarg{\epsilon_1,\ldots,\epsilon_d}{k_1,\ldots,k_d} & = \sum_{0 < n_1 < \cdots < n_d} \frac{(1 - (-1)^{n_1})\epsilon_1^{\frac{n_1-1}{2}}}{2 \, n_1^{k_1}} \cdots \frac{(1 - (-1)^{n_d})\epsilon_1^{\frac{n_d-1}{2}}}{2 \, n_d^{k_d}}  \\
		& = \frac{1}{2^d} \sum_{\substack{\eta_1^2 = \epsilon_1 \\\ldots \\\eta_d^2 = \epsilon_d}} \eta_1 \cdots \eta_d \zeta\sgnarg{\eta_1, \ldots, \eta_d}{k_1,\ldots,k_d}
	\end{align*}
	So alternating MtV's are expressed in terms of coloured MZV's of level \( N = 4 \); note however that \( \eta_1\cdots\eta_d = \pm 1 \) if an even number of \( \epsilon_i = -1 \), and \(\eta_1\cdots\eta_d = \pm \ii \) if an odd number of \( \epsilon_i = -1 \).  In particular, the MtV's correspond to the real part of purely real, respectively the imaginary part of purely imaginary, combinations in each case.
	\medskip
	
	Optimistically one expects some result of the following form. 
	
	\begin{Idea}\label{conj:tmm1m:basis}
		The stuffle-regularised alternating MtV's 
		\[
			\{ t^{\ast,V}(w) \mid w \in \{ 1, \overline{1} \}^\times \}
		\]		
		are linearly independent, and form a basis for the space of alternating MtV's.
	\end{Idea}

	Some issues do arise in attempting to investigate this, both from the motivic viewpoint, and from the classical viewpoint.  The issues should not be insurmountable, but we do postpone the motivic investigation for the moment.
	
	As just noted, for \( k_i \in \mathbb{Z}_{>0} \cup \overline{\mathbb{Z}_{>0}} \) (with \( \overline{k_i} \) denoting the argument \( k_i \) has associated sign \( \epsilon_i = -1 \)), one only has that 
	\[
		\ii^{\#\{ k_i \in \overline{\mathbb{Z}_{>0}} \}} t(k_1,\ldots,k_d) 
	\]
	is a \( \Q \)-linear combination of coloured MZV's of level \( N = 4 \).  So that one must extend the coefficients to \( \mathbb{Q}(\ii) \), in order to say anything about level \( N = 4 \) coloured MZV's from an MtV result, and vice-versa.
	
	Secondly, one must redevelop the background in \cite{murakami21}, and \cite{charltont2212} to extend the formulae for the motivic derivations \( D_{2r+1} \) to the case of alternating \( t \) values.  An deeper issue here is that there are more primitive elements at level \( N = 4 \), so that \( \ker D_{<N} = \zeta^\mot(N) \Q \oplus \zeta^\mot\sgnargsm{\ii}{N} \Q \) is no longer one-dimensional (see \cite[Corollary 5.1.3]{glanoisTh}).  Lifting identities by application of the period map requires identifying 2 unknown rational coefficients.  Working over \( \mathbb{Q} \), one can appeal to the real and imaginary parts for this, but after tensoring by \( \mathbb{Q}(\ii) \) this is not necessarily so straightforward.
	
	Finally, the natural regularisation \( V = \log(2) \), for \( t^{\ast,V}(1) \) immediately fails to give a basis in \autoref{conj:tmm1m:basis}.  One has that
	\begin{align*}
		t^{\ast,V}(1,1) & {} = \frac{1}{2} t^{\ast,V}(1)^2  \frac{1}{2} t(2) = \frac{1}{2} V^2 - \frac{\pi^2}{16} \\
		t(1,\overline{1}) & {} = \frac{1}{2} G - \frac{\pi}{8} \log(2) \\
		t^{\ast,V}(\overline{1},1) & {} = t(\overline{1})t^{\ast,V}(1) - t(1, \overline{1}) - t(\overline{2})  = -\frac{\pi}{4} V + \frac{1}{2} G + \frac{\pi}{8} \log(2) \\
		\quad\quad t(\overline{1},\overline{1}) & {} = \frac{1}{2} t(\overline{1})^2 - \frac{1}{2} t(2) = - \frac{\pi^2}{32}
	\end{align*}
	Here \( G = \beta(2) = 0.9159\ldots \) is the Catalan constant.  We refer to \cite[Section 6]{hoffman19} for these evaluations.  In the case \( V = \log(2) \), we have
	\[
		t^{\ast,V=\log(2)}(\overline{1}, 1) = -\frac{\pi}{8} \log(2) + \frac{1}{2} G  = t(1,\overline{1}) \,.
	\]
	So the stuffle-regularised \( t^{\ast,V}(w) \), \( w \in \{ 1,  \overline{1} \}^\times \) cannot be a basis for \( V =  \log(2) \).  The case \( V = \frac{1}{2} \log(2) \) does see to give a basis, though, and the regularisation \( t^{\ast,V=\frac{1}{2}\log(2)}(1) = \frac{1}{2} \log(2) \) is also a very natural one to take, based on the expression for \( t(n) = \frac{1}{2} \big( \zeta(n) - \zeta(\overline{n}) \big) \), extended to \( n = 1 \).  A more durable conjecture would be as follows.
	
	\begin{Conj}\label{conj:tmm1m:basis2}
		For any \( 0 < \lambda < 1 \), the stuffle-regularised alternating MtV's with \( V = \lambda \log(2) \), of the following form
		\[
		\{ t^{\ast,V}(w) \mid w \in \{ 1, \overline{1} \}^\times \}
		\]		
		are linearly independent, and form a basis for the space of alternating MtV's.
	\end{Conj}

	Those \( \lambda \) for which the regularisation \( t^{\ast,V} \), \( V = \lambda \log(2) \) does not give a basis in \autoref{conj:tmm1m:basis2}, should be termed \emph{singular regularisation parameters}.  The following \( V = \lambda \log(2) \) are singular regularisation parameters, first appearing at the indicated weight. \medskip
	\begin{center}
		\begin{tabular}{c|ccccccccccc}
			$\!\!N\!\!$ & 1 & 2 & 3 & 4 & 5 & 6 & 7 & 8 & 9 & 10 & 11 \\ \hline
			$ \!\!\lambda\!\!$ $\phantom{\mathllap{\displaystyle\frac{1}{1}}}$ &
			$0$ &
			$1$ &
			$-2$ &
			 $\frac{13}{11}$ &
			 \!${-}\frac{220}{203}$\! &
			 \!$\frac{4971}{3911}$\! &
			 \!${-}\frac{428854}{506177}$\! &
			 \!$\frac{8829285}{6699031}$\! &
			 \!${-}\frac{12070249400}{16117649299}$\! & 
			 \!$\frac{91040059801}{67506970721}$\! &
			 \!${-}\frac{917750647910294}{1321840200143647}\!$
			\end{tabular}
	\end{center}\medskip
	A new singular regularisation parameter \( \lambda \) appears in each weight, corresponding to a reduction of
	\[
 	t^{\ast,V}(\{\overline{1}\}^n,1) = \sum_{i=0}^{n-1} c_i t(\{\overline{1}\}^i, 1, \{\overline{1}\}^{n-i}) \,,
	\]
	for some \( c_i \in \Q \). 	For example, when \( V = \frac{13}{11} \log(2) \) in weight \( 4 \), we have
	\[
		t^{\ast, V=\frac{13}{11} \log(2)}(\overline{1}, \overline{1}, \overline{1}, 1) = \frac{41}{33} t(\overline{1}, \overline{1}, 1, \overline{1}) - \frac{15}{11} t(\overline{1}, 1, \overline{1}, \overline{1}) + \frac{15}{11} t(1, \overline{1}, \overline{1}, \overline{1}) \,,
	\]
	as can be verified with the evaluation in \autoref{thm:tmm1m}.  In weight \( N+1 \) this reduction can also be obtained directly from the identity in \autoref{thm:tmm1m}, when written in matrix form with rows indexed by \( c_i \) and columns by \( \zeta(2r+1) \) and \( \beta(2r) \); the factor \( (\frac{\pi}{4})^i \) is fixed by weight, so can be discarded.  The singular regularisation parameter in weight \( N+1 \) then corresponds to the determinant of the resulting matrix vanishing.
	
	The sequence \( (\lambda_i)_{i=1}^\infty \) of singular regularisation parameters appears to satisfy a number of properties.
	
	\begin{Conj}
		The sequence \( (\lambda_i)_{i=1}^\infty = (0, 1, -2, \frac{13}{11}, -\frac{220}{203}, \ldots) \) of singular regularisation parameters satisfies the following, discounting the \( \lambda_1 = 0 \) term:
		\begin{itemize}
			\item[i)] the sign of \( \lambda_i \) is \( (-1)^i \) for all \( i > 1 \),
			\item[ii)] the odd-indexed and even-indexed subsequences are increasing: \(\lambda_{2i+2}  > \lambda_{2i} \) and \( \lambda_{2i+3} > \lambda_{2i+1} \) for all \( i \geq 1 \).
			\item[iii)] the odd-indexed and even-indexed subsequences are bounded as follows: \( \lambda_{2i} \leq \frac{3}{2} \) and \( \lambda_{2i+1}\leq -\frac{1}{2} \) for all  \( i \geq 1 \).
			\item[iii$'$)] the odd-indexed and even-indexed subsequences have the following limits: \( \lim_{i\to\infty} \lambda_{2i} = \frac{3}{2} \) and \( \lim_{i\to\infty} \lambda_{2i+1} = -\frac{1}{2} \)
		\end{itemize}
	\end{Conj}

	As a final observation, based on the expression of the alternating MtV's as real/imaginary parts of coloured MZV's of level \( N = 4 \), depending on the parity of the number of \( \epsilon_i = -1 \), and the expected motivic results related to this, we should have the following.
	
	\begin{Conj}
		Any relation between alternating MtV's is homogeneous in the number of signs \( \epsilon_i = -1 \), counted modulo 2.
	\end{Conj}

	\begin{Eg}
		In weight 5, and with \( V = \frac{1}{2} \log(2) \), we have the identity
		\begin{align*}
			& t(\overline{1},2,\overline{2}) =  \\[1ex]
			&
			-\tfrac{143368}{3215}t(\overline{1},\overline{1},1,\overline{1},\overline{1})
			+\tfrac{121464}{3215}t(\overline{1},\overline{1},\overline{1},1,\overline{1})
			+\tfrac{183472}{16075}t^{\ast,V=\frac{1}{2}\log(2)}(\overline{1},\overline{1},\overline{1},\overline{1},1) \\
			& -48 t(1,\overline{1},\overline{1},\overline{1},\overline{1})			
			+48 t(\overline{1},1,\overline{1},\overline{1},\overline{1}) 
			 -\tfrac{24}{5} t^{\ast,V=\frac{1}{2}\log(2)}(1,\overline{1},\overline{1},1,1)
			+\tfrac{24}{5} t^{\ast,V=\frac{1}{2}\log(2)}(\overline{1},1,\overline{1},1,1) \\
			& - \tfrac{24}{25} t^{\ast,V=\frac{1}{2}\log(2)}(\overline{1},\overline{1},1,1,1)
			+ 4 t^{\ast,V=\frac{1}{2}\log(2)}(1,\overline{1},1,\overline{1},1)
			-4 t^{\ast,V=\frac{1}{2}\log(2)}(\overline{1},1,1,\overline{1},1) \,.
		\end{align*}
		One sees immediately that each term has either 2 or 4 barred entries, corresponding to 2 or 4 arguments with associated sign \( \epsilon_i = -1 \).
	\end{Eg}
	
	\bibliographystyle{habbrv}
	\bibliography{tmm1m}

\end{document}